\documentstyle{amsppt}
\magnification=1200
\hcorrection{.25in}
\TagsOnRight
\NoBlackBoxes
\topmatter
\NoBlackBoxes

\title Examples of Domains\\ 
with Non-Compact Automorphism Groups
\endtitle

\footnote[]{{\bf Mathematics
Subject Classification:} 32A07, 32H05, 32M05 \hfill}

\footnote[]{{\bf Keywords and Phrases:} Automorphism
groups, Reinhardt domains, circular domains.\hfill}

\author Siqi Fu, \ A. V. Isaev, \ and \ S. G. Krantz
\endauthor 

\abstract We give an example of a bounded, pseudoconvex, circular 
domain in ${\Bbb C}^n$ for any $n\ge 3$ with 
smooth real-analytic boundary and non-compact automorphism group, which 
is not biholomorphically equivalent
 to any Reinhardt domain. We also give an analogous example in ${\Bbb C}^2$, 
where the domain is bounded, non-pseudoconvex, is not
equivalent to any Reinhardt domain,  and the boundary is smooth   
real-analytic at all points except one.   
\endabstract
\endtopmatter  
\document
\def\qed{{\hfill{\vrule height7pt width7pt
depth0pt}\par\bigskip}}

Let $D$ be a bounded or, more generally, a hyperbolic
domain in ${\Bbb C}^n$. Denote by $\text{Aut}(D)$ the
group of biholomorphic self-mappings of $D$. The group
$\text{Aut}(D)$, with the topology given by uniform 
convergence on
compact subsets of $D$, is in fact a Lie group \cite{Kob}.

A domain $D$ is called Reinhardt if the standard action of 
the $n$-dimensional
torus ${\Bbb T}^n$ on ${\Bbb C}^n$,  
$$
z_j\mapsto e^{{i\phi}_j}z_j,\qquad {\phi}_j\in {\Bbb R},\quad
j=1,\dots,n,
$$
leaves $D$ invariant. For certain classes of domains with 
non-compact
 automorphism groups, Reinhardt domains serve as standard models up 
to biholomorphic equivalence (see e.g. \cite{R}, \cite{W}, 
\cite{BP}, \cite{GK1}, \cite{Kod}). 

It is an intriguing question whether {\it any} domain in 
${\Bbb C}^n$ with non-compact
 automorphism group and satisfying some natural geometric 
conditions is biholomorphically 
equivalent to a Reinhardt domain. The history of the study 
of domains with non-compact
 automorphism groups shows that there were expectations that 
the answer to this question
 would be positive (see \cite{Kra}). In this note we give 
examples that show that the answer 
is in fact negative.

While the domain that we shall consider in Theorem 1 
below has already been noted in the literature \cite{BP}, 
it has never been proved that 
this domain is not biholomorphically equivalent to a 
Reinhardt domain. 
Note that this domain is circular, i.e. 
it is invariant under the special rotations
$$
z_j\mapsto e^{{i\phi}}z_j,\qquad {\phi}\in {\Bbb R},\quad
j=1,\dots,n.
$$
Our first result is the following

\proclaim{Theorem 1} There exists a bounded, pseudoconvex, 
circular domain $\Omega\subset{\Bbb C}^3$
 with smooth real-analytic boundary and non-compact 
automorphism group, which is not 
biholomorphically equivalent to any Reinhardt domain.
\endproclaim
\demo{Proof} Consider the domain
$$
\Omega=\{|z_1|^2+|z_2|^4+|z_3|^4+(\overline{z_2}z_3+
\overline{z_3}z_2)^2<1\}.
$$
The domain $\Omega$ is invariant under the action of 
the two-dimensional torus ${\Bbb T}^2$
$$
\aligned
&z_1\mapsto e^{{i\phi}_1}z_1,\qquad {\phi}_1\in 
{\Bbb R},\\
&z_j\mapsto e^{{i\phi}_2}z_j,\qquad {\phi}_2\in 
{\Bbb R}, \quad j=2,3,
\endaligned 
$$
and therefore is circular. It is also a pseudoconvex, 
bounded domain with smooth real-analytic boundary. 
The automorphism group $\text{Aut}(\Omega)$  is 
non-compact since it contains the following subgroup
$$
\aligned
&z_1\mapsto\frac{z_1-a}{1-\overline{a}z_1},\\
&z_2\mapsto\frac{(1-|a|^2)^{\frac{1}{4}}z_2}
{(1-\overline{a}z_1)^{\frac{1}{2}}},\\
&z_3\mapsto\frac{(1-|a|^2)^{\frac{1}{4}}z_3}
{(1-\overline{a}z_1)^{\frac{1}{2}}},
\endaligned \tag{1}
$$
for a complex parameter $a$ with $|a|<1$.

We are now going to explicitly determine 
$\text{Aut}(\Omega)$. Let $F=(f_1,f_2,f_3)$ be an
 automorphism of $\Omega$. Then, since $\Omega$ is 
bounded, pseudoconvex and has real-analytic
 boundary, $F$ extends smoothly to 
$\overline{\Omega}$ \cite{BL}.
 Therefore, $F$ must preserve
 the rank of the Levi form ${\Cal L}_{\partial\Omega}(q)$ 
of $\partial\Omega$ at every $q\in\partial\Omega$.
 The only points where ${\Cal L}_{\partial\Omega}\equiv 0$  
are those of the form $(e^{i\alpha},0,0)$,
 $\alpha\in{\Bbb R}$.
These points must be preserved by $F$.
 This observation implies that $f_j(e^{i\alpha},0,0)=0$ for 
all $\alpha\in{\Bbb R}$, $j=2,3$.
 Restricting $f_2$, $f_3$ to the unit disc 
$\Omega\cap\{z_2=z_3=0\}$, we see that $f_j(z_1,0,0)=0$ for
 all $|z_1|\le 1$, $j=2,3$. Therefore, $F(0)=(b,0,0)$ for some $|b|<1$.
Taking the composition of $F$
 and the automorphism $G$ of the form (1) with $a=b$, 
we find that the mapping $G\circ F$ preserves the
 origin. Since $\Omega$ is circular, it follows from a theorem of 
H. Cartan \cite{C} that $G\circ F$ must be linear.
 Therefore, any automorphism of $\Omega$ is the composition of a linear automorphism and an automorphism of the form (1).

The above argument also shows that any linear 
automorphism of $\Omega$ can be written as
$$
\aligned
&z_1\mapsto e^{i\phi_1}z_1,\\
&z_2\mapsto az_2+bz_3,\\
&z_3\mapsto cz_3+dz_3,
\endaligned
$$
where $\phi_1\in{\Bbb R}$, $a,b,c,d\in {\Bbb C}$, 
and the transformation in the variables $(z_2, z_3)$ is
 an automorphism of the section $\Omega\cap \{z_1=0\}$. 
Further, since the only points of $\partial\Omega$
 where $\text{rank}\,{\Cal L}_{\partial\Omega}=1$  are 
those of the form $(z_1,w,\pm w)$ with $w\ne 0$ and since automorphisms of
 $\Omega$ preserve such points, it follows that any 
linear automorphism of $\Omega$ is in fact given by
$$
\aligned
&z_1\mapsto e^{i\phi_1}z_1,\\
&z_2\mapsto e^{i\phi_2}z_{\sigma(2)},\\
&z_3\mapsto \pm e^{i\phi_2}z_{\sigma(3)},
\endaligned
$$
where $\phi_1, \phi_2\in{\Bbb R}$, and $\sigma$ is a 
permutation of the set $\{2,3\}$.

The preceding description of $\text{Aut}(\Omega)$ 
implies that $\text{dim}\,\text{Aut}(\Omega)=4$.
That is to say, each of the four connected components of $\text{Aut}(\Omega)$  is parametrized by the point $a$ from the unit disc and by the 
rotation parameters $\phi_1, \phi_2$.

Suppose now that $\Omega$ is biholomorphically 
equivalent to a Reinhardt domain $D\subset{\Bbb C}^3$.
 Since $\Omega$ is bounded, it follows that $D$ is hyperbolic. 
It follows from \cite{Kru} that any hyperbolic Reinhardt
 domain $G\subset{\Bbb C}^n$ can be biholomorphically 
mapped onto its normilized form $\tilde G$ for which the identity
 component $\text{Aut}_0(\tilde G)$ of $\text{Aut}(\tilde G)$ is 
described as follows. There exist integers $0\le s\le
t\le p\le n$ and $n_i\ge 1$, $i=1,\dots,p$, with $\sum_{i=1}^p
n_i=n$, and real numbers ${\alpha}_i^k$, $i=1,\dots,s$,
$k=t+1,\dots,p$, and  ${\beta}_j^k$, $j=s+1,\dots,t$,
$k=t+1,\dots,p$  such that if we set
$z^i=\left(z_{n_1+\dots+n_{i-1}+1},\dots,
z_{n_1+\dots+n_i}\right)$,
$i=1,\dots,p$, then $\text{Aut}_0(\tilde G)$ is given by 
the mappings
$$
\aligned
&z^i\mapsto\frac{A^iz^i+b^i}{c^iz^i+d^i},
\quad i=1,\dots,s,\\
&z^j\mapsto B^jz^j+e^j,\quad j=s+1,\dots,t,\\
&z^k\mapsto
C^k\frac{\prod_{j=s+1}^t\exp\left(-\beta_j^k
\left(2\overline{e^j}^TB^jz^j+
|e^j|^2\right)\right)z^k}{\prod_{i=1}^s
(c^iz^i+d^i)^{2\alpha_i^k}},\quad
k=t+1,\dots,p,
\endaligned\tag{2}
$$
where
$$
\aligned
&\pmatrix
A^i&b^i\\
c^i&d^i
\endpmatrix\in SU(n_i,1),\quad i=1,\dots,s,\\
&B^j\in U(n_j),\quad e^j\in{\Bbb C}^{n_j},
\quad j=s+1,\dots,t,\\
&C^k\in U(n_k),\quad k=t+1,\dots,p.
\endaligned
$$
The normalized form $\tilde G$ is written as
$$
\aligned
G=\Biggl\{\left|z^1\right|&<1,\dots,\left|z^s\right|<1,\\
&\Biggl(\frac{z^{t+1}}{\prod_{i=1}^s
\left(1-\left|z^i\right|^2\right)^{{\alpha}_i^{t+1}}\prod_{j=s+1}^t
\exp\left(-{\beta}_j^{t+1}\left|z^j\right|^2\right)}\ ,\ \dots \ , \\
&\frac{z^{p}}{\prod_{i=1}^s
\left(1-\left|z^i\right|^2\right)^{{\alpha}_i^p}\prod_{j=s+1}^t
\exp\left(-{\beta}_j^p\left|z^j\right|^2\right)}\Biggr)\in\tilde
G_1\Biggr\},
\endaligned\tag{3}
$$
where $\tilde G_1:=\tilde G\bigcap\left\{z^i=0,\,
i=1,\dots,t\right\}$ is a hyperbolic Reinhardt domain in ${\Bbb
C}^{n_{t+1}}\times\dots\times{\Bbb C}^{n_p}$.

It is now easy to see that, for any hyperbolic 
Reinhardt domain $D\subset{\Bbb C}^3$ written in a 
normilized form $\tilde D$,
 $\text{Aut}_0(\tilde D)$ given by formulas (2) cannot 
have dimension equal to 4.

This completes the proof.\qed
\enddemo

{\bf Remark.} The theorem can be easily 
extended to ${\Bbb C}^n$ for any $n\ge 3$ 
(just replace $|z_1|^2$ in
 the defining function of $\Omega$ by 
$\sum_{j=1}^{n-2}|z_j|^2$, $z_2$ by $z_{n-1}$, $z_3$ by $z_n$). 

There is considerable evidence that, in complex dimension two, an example
such as that constructed in Theorem 1
does not exist.  Certainly the example provided above 
depends on the decoupling,
in the domain $\Omega$, of the variables $z_2, z_3$ from the
variable $z_1$.  Such decoupling is not possible when the dimension
is only two.

The work of 
Bedford and Pinchuk (see \cite{BP} and references therein)
suggests that the only 
smoothly bounded domains in ${\Bbb C}^2$ with non-compact automorphism groups
are (up to biholomorphic equivalence) the complex ellipsoids
$$
\Omega_{\alpha} = \{(z_1,z_2) \in 
{\Bbb C}^2: |z_1|^2 + |z_2|^{2\alpha} < 1\},
$$
where $\alpha$ is a positive integer. It is also a plausible 
conjecture that any bounded domain in ${\Bbb C}^2$ 
with non-compact automorphism group and a boundary of finite smoothness $C^k$ for $k\ge 1$, 
is biholomorphically equivalent to some $\Omega_{\alpha}$, where $\alpha\ge 1$ and is not necessarily an integer.   
Of course all the domains $\Omega_{\alpha}$ are pseudoconvex and Reinhard.

However, as the following theorem shows, if we allow the boundary to be non-smooth at just one point, then
the domain may be non-pseudoconvex and be non-equivalent to any Reinhardt domain.

\proclaim{Theorem 2} There exists a bounded, non-pseudoconvex domain $\Omega\subset{\Bbb C}^2$ 
with non-compact automorphism group such that $\partial\Omega$ 
is smooth real-analytic everywhere except one point (this
exceptional point is an orbit accumulation point for the automorphism group action), 
and such that $\Omega$ is not biholomorphically equivalent to any Reinhardt domain.
\endproclaim

For the proof of Theorem 2, we first need the following lemma.

\proclaim{Lemma A} If $\Omega\subset{\Bbb C}^2$ is a bounded, non-pseudoconvex, simply-connected domain such that the identity component $\text{Aut}_0(\Omega)$ of the automorphism group $\text{Aut}(\Omega)$ is non-compact, 
then $\Omega$ is not biholomorphically equivalent to any Reinhardt domain.
\endproclaim

\demo{Proof of Lemma A} Suppose that $\Omega$ is biholomorphically equivalent to a Reinhardt domain $D$. 
Since $\Omega$ is bounded, it follows that $D$ is hyperbolic. Also, since $\text{Aut}_0(\Omega)$ is 
non-compact, then so is $\text{Aut}_0(D)$. We are now going to show that any such domain $D$ is either 
pseudoconvex, or not simply-connected, or cannot be biholomorphically equivalent to a bounded domain. 
This result clearly implies the lemma.

We can now assume that the domain $D$ is written in its normalized form $\tilde D$ as in (3), and $\text{Aut}_0(\tilde D)$ is given by formulas (2). Then, since $\text{Aut}_0(\tilde D)$ is non-compact, it must be that $t>0$. Next, 
if $p=t$, then $\tilde D$ is either non-hyperbolic (for $s<t$), or (for $s=t$) is the unit ball or the unit polydisc and therefore is pseudoconvex. Thus we can assume that $t=1$, $p=2$, $n_1=n_2=1$.

Let $\tilde D_1\subset{\Bbb C}$ be the hyperbolic Reinhardt domain analogous to $\tilde G_1$ that was defined above (see (3)).
Clearly, there are the following possibilities for $\tilde D_1$:
\medskip 

\noindent {\bf (i)} $\tilde D_1=\{0<|z_2|<R\}$, $0<R<\infty$;

\noindent {\bf (ii)} $\tilde D_1=\{r<|z_2|<R\}$, $0<r<R\le \infty$;

\noindent {\bf (iii)} $\tilde D_1=\{|z_2|<R\}$, $0<R<\infty$.
\medskip 

For the cases {\bf (i)}, {\bf (ii)}, $\tilde D$ is always not simply-connected, and therefore we will concentrate on the case {\bf (iii)}.
 If $s=0$, then $\tilde D$ is not hyperbolic since it contains the complex line $\{z_2=0\}$.  
Thus we can assume that $s=1$. Next observe that, for $\alpha_1^2\ge 0$, 
the domain $\tilde D$ is always pseudoconvex.  Thus we may take $\alpha_1^2<0$. Then the domain $\tilde D$ has the form
$$
\tilde D=\left\{|z_1|<1,\,|z_2|<\frac{R}{(1-|z_1|^2)^{\gamma}}\right\},\qquad \gamma>0.
$$
We will now show that the above domain $\tilde D$ cannot be biholomorphically equivalent 
to a bounded domain. More precisely, we will show that any bounded holomorphic 
function on $\tilde D$ is independent of $z_2$.

Let $f(z_1,z_2)$ be holomorphic on $\tilde D$ and $|f|<M$ for some $M>0$. For every $\rho$ such that $|\rho|\le\frac{R}{2}$, the disc $\Delta_{\rho}=\{|z_1|<1,\,z_2=\rho\}$ is contained in $\tilde D$. We will show that ${\partial f}/\partial z_2 \equiv 0$  on every such $\Delta_{\rho}$, which implies that ${\partial f}/\partial z_2 \equiv 0$ everywhere in $\tilde D$.

Fix a point $(\mu,\rho)\in \Delta_{\rho}$ and restrict $f$ to the disc $\Delta'_{\mu}=\{z_1=\mu,\,|z_2|<R_{\mu}\}$, where $R_{\mu}={R}/{2(1-|\mu|^2)^{\gamma}}$ . Clearly, $(\mu,\rho)\in\Delta'_{\mu}$ and $\overline{\Delta'_{\mu}}\subset \tilde D$. By the Cauchy Integral Formula
$$
f(\mu,z_2)=\frac{1}{2\pi i}\int_{\partial \Delta'_{\mu}}\frac{f(\mu,\zeta)}{\zeta-z_2}\,d\zeta,
$$
for $|z_2|<R_{\mu}$, and therefore
$$
\frac{\partial f}{\partial z_2}(\mu,\rho)=\frac{1}{2\pi i}\int_{\partial \Delta'_{\mu}}\frac{f(\mu,\zeta)}{(\zeta-\rho)^2}\,d\zeta.
$$
Hence
$$
\left|\frac{\partial f}{\partial z_2}(\mu,\rho)\right|\le\frac{MR_{\mu}}{(R_{\mu}-|\rho|)^2}.
$$
Letting $|\mu|\rightarrow 1$ and taking into account that $R_{\mu}\rightarrow\infty$, we see that
$\left|{\partial f}/\partial z_2 (\mu,\rho)\right|\rightarrow 0$ as  $|\mu|\rightarrow 1$. Therefore, ${\partial f}/\partial z_2 \equiv 0$ on $\Delta_{\rho}$.

The lemma is proved.\qed
\enddemo

\demo{Proof of Theorem 2} We will now present a domain that satisfies the conditions of the lemma. Set
$$
\Omega=\left\{|z_1|^2+|z_2|^4+8|z_1-1|^2\left(\frac{z_2^2}{z_1-1}-\frac{3}{2}\frac{|z_2|^2}{|z_1-1|}+\frac{\overline{z_2}^2}{\overline{z_1}-1}\right)^2<1\right\}.
$$

The domain $\Omega$ is plainly bounded since the third term on the left is non-negative. Next, the identity component $\text{Aut}_0(\Omega)$ of its automorphism group is non-compact since it contains the subgroup
$$
\aligned
&z_1\mapsto\frac{z_1-a}{1-az_1},\\
&z_2\mapsto\frac{(1-a^2)^{\frac{1}{4}}z_2}
{(1-az_1)^{\frac{1}{2}}},
\endaligned
$$
where $a\in(-1,1)$. 

Further, $\Omega$ is simply-connected, since the family of mappings $F_{\tau}(z_1,z_2)=(z_1,\tau z_2)$, $0\le\tau\le 1$, 
retracts $\Omega$ inside itself, as $\tau\rightarrow 0$, to the unit disc $\{|z_1|<1,\, z_2=0\}$ (which is simply-connected).

To show that $\Omega$ is not pseudoconvex, consider its unbounded realization. Namely, under the mapping
$$
\aligned
&z_1\mapsto \frac{z_1+1}{z_1-1},\\
&z_2\mapsto \frac{\sqrt{2}z_2}{\sqrt{z_1-1}},
\endaligned\tag{4}
$$
the domain $\Omega$ is transformed into the domain
$$
\Omega'=\left\{\text{Re}\,z_1+\frac{1}{4}|z_2|^4+2\left(z_2^2-\frac{3}{2}|z_2|^2+\overline{z_2}^2\right)^2<0\right\}.
$$
It is easy to see that at the boundary point $(-\frac{3}{4},1)\in\partial\Omega'$ the Levi form of $\partial\Omega'$ is equal to $-|z_2|^2$, and thus is negative-definite. Therefore, $\Omega$ is non-pseudoconvex.

Hence, by Lemma A, $\Omega$ is not biholomorphically equivalent to any Reinhardt domain. 

Next, if $\phi$ denotes the defining function of $\Omega$, the following holds at every boundary point 
of $\Omega$ except $(1,0)$:
$$
\frac{\partial \phi}{\partial z_1}=\frac{1}{z_1-1}\left(-\frac{z_2}{2}\frac{\partial \phi}{\partial z_2}+1-\overline{z_1}\right),
$$
and therefore $\text{grad}\,\phi$ does not vanish at every such point. Hence, $\partial\Omega$ is smooth real-analytic everywhere except at $(1,0)$.

The theorem is proved.\qed
\enddemo

{\bf Remarks.}   \break
\smallskip

\noindent {\bf 1.}  The hypothesis of 
simple connectivity in Lemma A is automatically satisfied if, for example, the boundary of the domain is locally variety-free and smooth near some orbit accumulation point for the automorphism group of the domain (see e.g. \cite{GK2}). For a smoothly bounded domain it would follow from a conjecture of Greene/Krantz \cite{GK3}.

\noindent {\bf 2.}  Tedious calculations show that the 
boundary of the domain $\Omega$ in Theorem 2 is 
quite pathological near the exceptional point $(1,0)$.  
It is not Lipschitz-smooth of any positive degree.  It would be
interesting to know whether there is an example with Lipschitz-1
boundary at the bad point. 

In fact, many more examples similar to that in Theorem 2 can be
constructed in the following way.  Let 
$$
\Omega' = \{(z_1,z_2) \in {\Bbb C}^2: 
\hbox{Re}\, z_1 + P(z_2) < 0\},\tag{5}
$$ 
where $P=|z_2|^{2m}+Q(z_2)$ is a homogeneous 
non-plurisubharmonic polynomial, $m$ is a positive 
integer, and $Q(z_2)$ is positive away from the origin. 
Then, by a mapping analogous to (4), $\Omega'$ can be transformed
into a bounded domain $\Omega$.  The domain $\Omega$ is simply-connected,
non-pseudoconvex, $\text{Aut}_0(\Omega)$ is non-compact, and
$\partial \Omega$ is smooth real-analytic everywhere except at the
point $(1,0)$.  For all such examples, $\partial \Omega$
is not Lipschitz-smooth of any positive degree at $(1,0)$.

It is also worth noting that, in the example contained in Theorem 2, the
point $(-1,0)$ is also an orbit accumulation point, but $\partial \Omega$
is smooth real-analytic at this point. 

\noindent {\bf 3.} It is conceivable that
the domain $\Omega$ as in Theorem 2 has an alternative, smoothly bounded
realization, but it looks plausible that if in formula (5) we allow $P(z_2)$ to be an arbitrary homogeneous polynomial positive away from the origin with no harmonic term, then domain (5) does not have a bounded realization with $C^1$-smooth boundary, unless $P(z_2)=c|z_2|^{2m}$, where $c>0$ and $m$ is a positive integer. 
\bigskip

{\bf Acknowledgements.} This work was completed while the second 
author was an Alexander von Humboldt Fellow at the 
Research at MSRI by the third author supported in part by 
NSF Grant DMS-9022140.
University of Wuppertal.

\bigskip
\pagebreak

{\obeylines
\font\foo=cmcsc9
\foo
Siqi Fu
Department of Mathematics 
University of California, Irvine, CA 92717
USA 
E-mail address: sfu\@math.uci.edu
\bigskip 
A. V. Isaev 

Centre for Mathematics and Its Applications 
The Australian National University 
Canberra, ACT 0200
AUSTRALIA 
E-mail address: Alexander.Isaev\@anu.edu.au
\smallskip
and
\smallskip
Bergische Universit\"at
Gesamthochschule Wuppertal
Mathematik (FB 07)
Gaussstrasse 20
42097 Wuppertal
GERMANY
E-mail address: Alexander.Isaev\@math.uni-wuppertal.de
\bigskip
S. G. Krantz
Department of Mathematics
Washington University, St.Louis, MO 63130
USA 
E-mail address: sk\@math.wustl.edu
\smallskip
and
\smallskip
MSRI
1000 Centennial Drive
Berkeley, California 94720
USA
E-mail address: krantz\@msri.org

\enddocument